\documentclass{article}
\usepackage{amsmath}
\usepackage{amsfonts}
\usepackage{amsthm}
\usepackage{centernot}
\usepackage{url}

\def\intZ{ \mathbb{Z} }

\def\natN{  \mathbb{N} }

\def\mod{\bmod}

\newtheorem{Conjecture}{Conjecture}[section]

\def\Jacobi#1#2{{\left( \frac{#1}{#2} \right)}}

\title{A Survey on the Ternary Purely Exponential Diophantine Equation $a^x + b^y = c^z$  }
\author{Maohua Le \and Reese Scott \and Robert Styer}
\date{7 November 2018}

\begin{document}

\maketitle

%\begin{abstract}
Let $a$, $b$, $c$ be fixed coprime positive integers with $\min\{a,b,c\}>1$.  In this survey, we consider some unsolved problems and related works concerning the positive integer solutions $(x,y,z)$ of the ternary purely exponential diophantine equation $a^x + b^y = c^z$.  
%\end{abstract} 

\bigskip



\section{Introduction}  %1

Let $\intZ$, $\natN$ be the sets of all integers and positive integers, respectively. A great deal of number theory arises from the discussion of the integer or rational solutions of a polynomial equation with integer coefficients.  Such equations are called diophantine equations.  Using the definition given by T. N. Shorey and R. Tijdeman \cite{ShTi1}, those diophantine equations with fixed bases and variable exponents are called purely exponential equations.  Let $a$, $b$, $c$ be fixed coprime positive integers with $\min\{a,b,c\}>1$.  In this survey we investigate the ternary purely exponential equations of the form 
$$ a^x + b^y = c^z, x,y,z \in \natN. \eqno{(1.1)} $$

We give an outline of the contents of this survey.  

First, in this Introduction (Section 1), after some explanatory comments to help orient the reader, we give the current bounds on the variable exponents $x$, $y$, and $z$.  We then give the current bounds on the number of solutions $(x,y,z)$ to (1.1).  

Section 2 deals with the Je\'smanowicz conjecture, which states that, if $(a,b,c)$ is a Pythagorean triple, then the only solution to (1.1) is $(x,y,z) = (2,2,2)$.  Work on this conjecture is characterized by many different highly specific results, making it difficult for those working on this conjecture to know if they are making full use of what is available and not reproving old results. Here the first author presents the existing work in an organized manner, which we hope will greatly facilitate future work on this problem.    

Section 3 introduces the Terai-Je\'smanowicz conjecture, which states that (1.1) has at most one solution with $\min\{x,y,z\}>1$. Here the first author clarifies the various oversights which have occurred in previous formulations of the conjecture.  Then a third conjecture, which includes both the Je\'smanowicz and Terai-Je\'smanowicz conjectures, is stated and briefly discussed.  The remainder of Section 3 is largely concerned with various cases for which it can be shown that (1.1) has at most one solution.   

Before proceeding, it might be helpful to point out what is {\it not} included in this survey.  

First, when one of $x$ or $y$ is fixed, (1.1) becomes the familiar Pillai equation, which has already been thoroughly handled in an excellent survey by M. Waldschmidt \cite{Wald}.  So here we treat only cases in which all of $x$, $y$, and $z$ are variable.  Second, (1.1) can be viewed as a specific case of an $S$-unit equation.  The many results of this type are usually not included here.  Third, we do not reference any work from the vast and familiar literature dealing with cases in which at least one of the bases $a$, $b$, $c$ is variable.  Finally, we deal largely with fairly recent results; we do not attempt to give a complete history of the problem.  

We note that the first author conceived the project and drafted the bulk of this paper and its extensive bibliography, the second author contributed some sections showing connections to more general work on exponential equations as well as contributing multiple improvements, and the third author edited previous versions and prepared the final version as well as handling various technical matters.  

We now begin our discussion of (1.1).  

In 1933, K. Mahler \cite{Mah1} used his $p$-adic analogue of the diophantine approximation method of Thue-Siegel to prove that (1.1) has only finitely many solutions $(x,y,z)$.  His method is ineffective.  An effective result for solutions of (1.1) was given by A. O. Gel'fond \cite{Gel1}.  In 1999, M.-H. Le \cite{Le7} used some elementary methods to prove that if $c \equiv 1 \bmod 2$, then the solutions $(x,y,z)$ of (1.1) satisfy $z < \frac{2}{\pi} a b \log(2 e a b)$.  This gives a bound on $z$ independent of $c$.  Recently, R. Scott and R. Styer \cite{ScSt4} proved that if $c \equiv 1 \bmod 2$ then $z < \frac{1}{2} ab$.  In some cases, finding a bound on $z$ independent of $c$ can be taken quite a bit further: when every prime dividing $ab$ is in a given finite set of primes $S$, it is often possible to show that (1.1) implies $z=1$ except for a finite list of specified exceptions; to do this, elementary methods sometimes suffice (e. g., $S= \{ 2,3,5 \}$), but in general they do not, particularly when $2 \not\in S$.  M. A. Bennett and N. Billerey \cite{BeBi} deal with the case $S = \{ 3,5,7 \}$, completely handling not only (1.1) but also the more general case in which $A$ and $B$ are $S$-units such that $A+B=c^z$, $z>1$.   

We now consider finding bounds on $\max\{x,y,z\}$.  N. Hirata-Kohno \cite{Hir1} showed that if $c \equiv 1 \bmod 2$ then $\max\{x,y,z\} < 2^{288} \sqrt{abc} (\log(abc))^3$.   For general $a$, $b$, and $c$, combining a lower bound for linear forms in two logarithms and an upper bound for the $p$-adic logarithms, Y.-Z. Hu and M.-H. Le \cite{HuLe1} proved that 
$$ \max\{x,y,z\} < 155000 (\log(\max\{a,b,c\}))^3.  \eqno{(1.2)} $$
Very recently they improved that to 
 $$ \max\{x,y,z\} < 6500 (\log(\max\{a,b,c\}))^3.  \eqno{(1.3)} $$
(see Y.-Z. Hu and M.-H. Le \cite{HuLe2}).  

Let $N(a,b,c)$ denote the number of solutions $(x,y,z)$ of (1.1).  As a straightforward consequence of an upper bound for the number of solutions of binary $S$-unit equations due to F. Beukers and H. P. Schlickewei \cite{BeSc1}, we have $N(a,b,c)\le 2^{36}$.  Considering the parity of $x$ and $y$, all solutions $(x,y,z)$ of (1.1) can be put into the following four classes:

Class I:  $x \equiv y \equiv 0 \bmod 2$.

Class II: $x \equiv 1 \bmod 2$ and $y \equiv 0 \bmod 2$.

Class III: $x \equiv 0 \bmod 2$ and $y \equiv 1 \bmod 2$.

Class IV:  $x \equiv y \equiv 1 \bmod 2$.  

Under this approach, M.-H. Le \cite{Le7} proved that if $a<b$ and $c \equiv 1 \bmod 2$, then each class has at most $2^{\omega(c) - 1}$ solutions $(x,y,z)$, where $\omega(c)$ is the number of distinct prime divisors of $c$, except for $(a,b,c) = (3, 10, 13)$.  This implies that $N(a,b,c) \le 2^{\omega(c)+1} $ if $c \equiv 1 \bmod 2$.  Recently, R. Scott and R. Styer \cite{ScSt4} improved this: if $c \equiv 1 \bmod 2$, all solutions $(x,y,z)$ of (1.1) occur in at most two distinct parity classes and each class has at most one solution when solutions exist in two distinct classes, and at most two solutions otherwise.  Therefore, we have $N(a,b,c) \le 2$ if $c \equiv 1 \bmod 2$.  Very recently, combining the upper bound (1.3) with some elementary methods, Y.-Z. Hu and M.-H. Le \cite{HuLe3, HuLe2} proved that if $\max\{a,b,c\} > 5 \times 10^{27}$ then $N(a,b,c) \le 3$.  Moreover, if $2 \mid c$ and $\max\{ a,b,c \} \ge 10^{62}$, then $N(a,b,c) \le 2$.  

Notice that, for any positive integer $k$ with $k \ge 2$, if $(a,b,c) = (2, 2^k -1, 2^k + 1)$, then (1.1) has two solutions $(x,y,z) = (1,1,1)$ and $(k+2,2,2)$.  Hence, as observed by R. Scott and R. Styer \cite{ScSt4}, there exist infinitely many triples $(a,b,c)$ with $N(a,b,c) =2$.  
Therefore, in general, $N(a,b,c) \le 2$ may be the best upper bound for $N(a,b,c)$.  

By combining the results in \cite{HuLe2} and \cite{ScSt4}, one reduces the problem of finding $N(a,b,c)>2$ to a finite search.  

Because computer searches suggest that $N(a,b,c) \le 1$ except for the known cases of double solutions, there have been a series of conjectures concerning exact upper bounds for $N(a,b,c)$.  In the next two sections, we shall introduce these conjectures and related works.

\section{Je\'smanowicz conjecture}  %2 

A positive integer triple $(A,B,C)$ is called a Pythagorean triple if 
$$ A^2 + B^2 = C^2.  \eqno{(2.1)} $$
In particular, if $\gcd(A,B) = 1$, then $(A,B,C)$ is called a primitive Pythagorean triple.  By (2.1), every Pythagorean triple $(A,B,C)$ can be expressed as 
$$ A = A_1 n, B = B_1 n, C = C_1 n, n \in \natN, \eqno{(2.2)} $$
where $(A_1, B_1, C_1)$ is a primitive Pythagorean triple.  Notice that every primitive Pythagorean triple $(A_1, \allowbreak B_1, \allowbreak  C_1)$ satisfies $C_1 \equiv 1 \bmod 2$, also $A_1$ and $B_1$ have opposite parity.  We may therefore assume without loss of generality that $A_1 \equiv 1 \bmod 2$ and $B_1 \equiv 0 \bmod 2$.  So we have 
$$ A_1 = f^2 - g^2, B_1 = 2 f g, C_1 = f^2 + g^2, f,g \in \natN, f>g, \gcd(f,g) = 1, fg \equiv 0 \bmod 2 \eqno{(2.3)} $$
(see L. J. Mordell \cite[Chapter 3]{Mor1}).  

In 1956, L. Je\'smanowicz \cite{Jes1} conjectured that if $(A,B,C)$ is a Pythagorean triple, then the equation 
$$ A^x + B^y = C^z, x,y,z \in \natN \eqno{(2.4)} $$
has only one solution $(x,y,z) = (2,2,2)$.  By (2.2), (2.3), and (2.4), Je\'smanowicz' conjecture can be rewritten as follows:

\begin{Conjecture}  %2.1
For any positive integer $n$, the equation 
$$ ((f^2 - g^2)n)^x + (2fgn)^y = ((f^2 + g^2)n)^z, x,y,z \in \natN \eqno{(2.5)} $$
has only the single solution $(x,y,z) = (2,2,2)$, where $f$, $g$ are positive integers satisfying (2.3).  
\end{Conjecture}

Je\'smanowicz' conjecture has been proven to be true in many special cases.  But, in general, the problem is not solved as of yet.

\subsection{Primitive cases}  % 2.1 

In this subsection we consider Conjecture 2.1 for $n=1$.  Many of the early works on Conjecture 2.1 deal with (2.5) for the case 
$$ n = 1, f = g+1.  \eqno{(2.6)} $$
When (2.6) holds, Conjecture 2.1 is true if one of the following conditions is satisfied:

(i) (W. Sierpi\'nski \cite{Sie1}) $g=1$.

(ii)  (L. Je\'smanowicz \cite{Jes1}) $2 \le g \le 5$.  

(iii)  (Z. Ke \cite{Ke1}, \cite{Ke2}) $g \equiv 1,3,4,5,7,9,10,11 \mod 12$, $g \equiv 2 \bmod 5$, $g \equiv 3 \bmod 7$, $g \equiv 4 \bmod 9$, $g \equiv 5 \bmod 11$, $g \equiv 6 \bmod 13$, or $ g \equiv 7 \mod 15$.  

(iv)  (D.-M. Rao \cite{Rao1}) $g \equiv 2, 6 \bmod 12$.

(v)  (Z. Ke \cite{Ke3}, \cite{Ke4}, \cite{Ke5}); Z. Ke and Q. Sun \cite{KeSu1}) $g \le 6144$.  

In 1965, V. A. Dem'janenko \cite{Dem1} completely solved the case (2.6).  Four decades later, Y.-Z. Hu and P.-Z. Yuan \cite{HuYu3} gave a new proof of Dem'janenko's result.  Moreover, Z. Ke, T. J\'ozefiak, J.-R. Chen, V. D. Podyspanin, Z.-F. Cao, W.-J. Chen, A. Grytczuk and A. Grelak successively discussed some rather special cases of Conjecture 2.1 for $n=1$.  The proofs of the above results are elementary.  A detailed record can be found in Z. Ke and Q. Sun \cite[Section 7.1]{KeSu2}, Z.-F. Cao \cite[Section 9.2]{Cao6}, G. Soydan, M. Demirci, I. N. Cangul, and A. Togb\'e \cite{SDCT1}.  

For any positive integer $t$, let $P(t)$ denote the product of the distinct prime divisors of $t$.  

Except for the above mentioned works, the existing results on Conjecture 2.1 for $n=1$ can be divided into three cases as follows. 

\smallskip

Case I.  A solution $(x,y,z)$ of (2.5) with $(x,y,z) \ne (2,2,2)$ is called exceptional.  If $n=1$, then the exceptional solutions of (2.5) have the following properties:

(i)  (Z. Li \cite{ZLi1}) If $x \equiv y \equiv z \equiv 0 \bmod 2$, then $x \equiv y \equiv z \equiv 2 \bmod 4$. 

(ii)  (T. Miyazaki \cite{Miy4}) If $y \equiv z \equiv 0 \bmod 2$, then $y \equiv z \equiv 2 \bmod 4$.  

(iii)  (M.-J. Deng and D.-M. Huang \cite{DeHu1}) If $f \equiv 2 \bmod 4$, $g \equiv 3 \bmod 4$, and $f+g \equiv 1 \bmod 16$, then $y=1$.  

(iv)  (M.-M. Ma and Y.-G. Chen \cite{MaCh1})  If $fg \equiv 2 \bmod 4$, then $y=1$.  

\smallskip

Case II.  If $f$ and $g$ take the following values, then Conjecture 2.1 is true for $n=1$.  

(i)  (W.-D. Lu \cite{Lu1}) $g=1$.

(ii)  (H.-L. Liu \cite{HLLiu1}); X.-F. An \cite{XFAn1}) $g=2$ and $f$ is an odd prime power.  

(iii)  (N. Terai \cite{Ter8})   $g=2$.  

(iv)  (C.-Y. Fu \cite{Fu1})  $g=6$ and $f$ is an odd prime power.  

(v)  (T. Miyazaki \cite{Miy7})  $f^2 - 2fg - g^2 = \pm 1$.

\smallskip 

Case III.  If $f$ and $g$ satisfy the following congruence and divisibility conditions, then Conjecture 2.1 is true for $n=1$.  

(i)  (K. Takakuwa \cite{Tak1})  $f \equiv 2 \bmod 4$ and $g \in \{ 3,7,11,15 \}$.

(ii)  (M.-H. Le \cite{Le3})  $fg \equiv 2 \bmod 4$ and $f^2 + g^2$ is an odd prime power.  

(iii)  (Y.-D. Guo and M.-H. Le \cite {GuLe1}; M.-H. Le \cite{Le4})  $f \equiv 2 \bmod 4$, $g \equiv 3 \bmod 4$, and $f > 81 g$.  

(iv)  (K Takakuwa and Y. Asaeda \cite{TaAs1}, \cite{TaAs2}, \cite{TaAs3}) $f \equiv 2 \bmod 4$, $g$ is an odd prime with $g \equiv 3 \bmod 4$, and the divisors of $f^2 - g^2$ satisfy certain conditions.  

(v)  (M.-J. Deng \cite{Den1}, \cite{Den2}; J.-H. Wang and M.-J. Deng \cite{WaDe1}; M.-J. Deng and G. L. Cohen \cite{DeCo2}; T. Miyazaki \cite{Miy1}; S.-Z. Li \cite{SZLi1}; Y. An \cite{An1}; J.-J. Xing \cite{Xin1}; S.-S. Gou \cite{Gou1}; C.-Y. Zheng \cite{Zhe1}) The divisors of $2fg$ and $f^2+g^2$ satisfy certain congruence properties.  

(vi)  (W.-J. Guan \cite{WJGua1}) $(f,g) = (1,6)$, $(2,5)$, $(5,2)$, $(6,1) \bmod 8$.

(vii)  (T. Miyazaki \cite{Miy3})  $f^2 - g^2 \equiv \pm 1 \bmod 2fg$ or $f^2 + g^2 \equiv 1 \bmod 2fg$.  

(viii)  (T. Miyazaki, P.-Z. Yuan, and D.-Y. Wu \cite{MYW1})  $f^2 - g^2 \equiv \pm 1 \bmod P(2fg)$ or $f^2 + g^2 \equiv \pm 1 \bmod P(2fg)$.  

(ix)  (Y. Fujita and T. Miyazaki \cite{FuMi1})  $2fg \equiv 0 \bmod 2^k$ and $2fg \equiv \pm 2^k \bmod (f^2 - g^2)$, where $k$ is a positive integer with $k \ge 2$.  

(x)  (Y. Fujita and T. Miyazaki \cite{FuMi2})  $2fg$ has a divisor $d$ with $d \equiv \pm 1 \bmod (f^2 - g^2)$.  

(xi)  (Q. Han and P.-Z. Yuan \cite{HaYu1})  $fg \equiv 2 \bmod 4$ and $f+g$ has a prime divisor $p$ with $p \not\equiv 1 \bmod 16$.  

(xii)  (M.-J. Deng and J. Guo \cite{DeGu1})  $g \equiv 2 \bmod 4$ and $n<600$.  

%  evidently this is not correct  (xv)  (Y.-F. He \cite{He1})  $fg \equiv 2 \bmod 4$.  

(xiii)  (T. Miyazaki and N. Terai \cite{MiTe2})  $f> 72 g$, $g \equiv 2 \bmod 4$, and $g/2$ is an odd prime power or a square.  

(xiv)  (P.-Z. Yuan and Q. Han \cite{YuHa1}) $fg \equiv 2 \bmod 4$,  $f>72 g$, and the divisors of $f$, $g$ satisfy some conditions.

(xv)  (M.-H. Le \cite{Le19})  $fg \equiv 2 \bmod 4$ and $f > 30.8 g$.     

(xvi)  (M.-H. Le \cite{Le17})    $f^2+g^2 > 4 \times 10^9$ and $\gcd(f^2+g^2, ((f^2-g^2)^l -\lambda) / (f^2+g^2) ) = 1$, where $l$ is the least positive integer with $(f^2-g^2)^l \equiv \lambda \bmod (f^2+g^2)$, $\lambda \in \{1, -1\}$.

\subsection{Non-primitive cases}  % 2.2

In this subsection we consider Conjecture 2.1 for $n>1$.  There was no study of this problem until 1998.  In 1998, M.-J. Deng and G. L. Cohen \cite{DeCo1} proved that if $n>1$, $f^2-g^2$ is an odd prime power, and either $n\equiv 0 \bmod P(2fg)$ or $2fg \not\equiv 0 \bmod P(n)$, then (2.5) has only one solution $(x,y,z) = (2,2,2)$.  One year later, M.-H. Le \cite{Le6} gave a more general result on Conjecture 2.1 for $n>1$.  He proved that if $n>1$, then every exceptional solution $(x,y,z)$ of (2.5) satisfies one of the following conditions:  

(i)  $\max\{x,y\} > \min\{x,y \} > z$, $f^2 +g^2 \equiv 0 \bmod P(n)$, and $P(n) < P(f^2+g^2)$.  

(ii)  $x>z>y$ and $2fg \equiv 0 \bmod P(n)$.

(iii)  $y>z>x$ and $f^2-g^2 \equiv 0 \bmod P(n)$.  

Sixteen years later, H. Yang and R.-Q. Fu \cite{YaFu2} simplified Le's result, showing that the condition (i) can be removed.  
The above two results are invariably used for solving (2.5) with $n>1$.  

Many works investigated (2.5) in the case
$$n>1, f=2^k, g=1, k \in \natN.  \eqno{(2.7)} $$ 
When (2.7) holds, if $k$ takes the following values, then Conjecture 2.1 is true.

(i)  (M.-J. Deng and G. L. Cohen \cite{DeCo1}) $k=1$.

(ii)  (Z.-J. Yang and M. Tang \cite{YaTa1}) $k=2$.

(iii)  (S.-S. Gou and H. Zhang \cite{GoZh1})  $k=3$.  

(iv)  (M. Tang and Z.-J. Yang \cite{TaZJYa1}) $k \in \{1,2,4,8\}$.

(v)  (M.-J. Deng \cite{Den4})  $ k \in \{2,3,4,5\}$.

(vi)  (M. Tang and J.-X. Weng \cite{TaWe1})  $k = 2^s$, where $s$ is a nonnegative integer.  

In about 2014, X.-W. Zhang and W.-P. Zhang \cite{ZhZh1}, and T. Miyazaki \cite{Miy10} independently solved the case (2.7), namely, they proved that if $n$, $f$, and $g$ satisfy (2.7), then Conjecture 2.1 is true.  

Conjecture 2.1 for $n>1$ has also been verified in the following cases.  

(i)  (Z.-J. Yang and J.-X. Weng \cite{YaWe1})  $(f,g) = (6,1)$. 

(ii)  (G. Soydan, M. Demirci, I. N. Cangul, and A. Togb\'e \cite{SDCT1})  $(f,g) = (10, 1)$.

(iii)  (B.-L. Liu \cite{BLLiu3})  $(f,g)=(12,1)$.  

(iv)  (D.-R. Ling and J.-X. Weng \cite{LiWe1})  $(f,g) = (14,1)$.  

(v)  (M.-M. Ma and J.-D. Wu \cite{MaWu2}; C.-F. Sun and Z. Cheng \cite{SuCh3}) $g=1$, either $n \equiv 0 \bmod P(f^2-1)$ or $f=2 p^k$ and $4 p^{2k} - 1 \not\equiv 0 \bmod P(n)$, where $p$ is an odd prime with $p \equiv 3 \bmod 4$, $k$ a positive integer.  

(vi)  (Z. Cheng, C.-F. Sun, and X.-N. Du \cite{CSD1})  $(f,g) = (5,2)$. 

(vii)  (C.-F. Sun and Z. Cheng \cite{SuCh1}; G. Tang \cite{Tan1})  $(f,g) = (7,2)$. 

(viii)  (C.-Y. Fu and M.-J. Deng \cite{FuDe1})  $(f,g) = (7,2)$ or $(49,2)$.  

(ix)  (C.-F. Sun and Z. Cheng \cite{SuCh2})  $(f,g) = (9,2)$.  

(x)  (W.-Y. Lu, L. Gao, and H.-F. Hao \cite{LGH1})  $(f,g) = (11,2)$.  

(xi)  (H. Yang, R.-Z. Ren, and R.-Q. Fu \cite{YRF1}; F.-J. Chen \cite{FJChe1})  If $f$ is an odd prime and $g=2$, then (2.5) has no solutions $(x,y,z)$ with $x>z>y$.  

(xii)  (Y.-M. Li \cite{YMLi1})  $(f,g) = (8,3)$.  

(xiii)  (M.-M. Ma \cite{Ma1}; M.-M. Ma and J.-D. Wu \cite{MaWu1})  $(f,g) = (9,4)$.  

(xiv)  (M.-J. Deng and G. L. Cohen \cite{DeCo1})  $(f,g) = (3,2)$, (4,3)$, (5,4)$, $(6,5)$.  

(xv)  (M.-J. Deng \cite{Den3})  $(f,g) = (7,6)$, $(8,7)$.

(xvi)  (H. Che \cite{Che1})  $(f,g) = (11,10)$.  

(xvii)  (H.-N. Sun \cite{Sun1})  $(f,g)=(18,17)$.  

(xviii)  (L.-L. Wang \cite{Wan1})  $(f,g)=(20,19)$. 

(xix)  (C.-N. Lin \cite{CNLin1})  $(f,g) = (26,25)$.

(xx)  (J. Ma \cite{JMa1})  $(f,g) = (29,28)$.  

(xxi)  (W.-Y. Lu, L. Gao, X.-H. Wang, and H.-F. Hao \cite{LGWH1})  $(f,g) = (46,45)$.  

(xxii)  (H. Yang and R.-Q. Fu \cite{YaFu2})  $(f,g) = (2^k, 2^k-1)$ and $2^k-1$ is an odd prime, where $k$ is a positive integer.  

(xxiii)  (H. Yang and R.-Q. Fu \cite{YaFu3})  $(f,g) = (2^k + 1, 2^k)$ and $2^k+1$ is an odd prime, where $k$ is a positive integer.  

(xxiv)  (T.-T. Wang, X-H. Wang, and Y.-Z. Jiang \cite{WWJ1})  If $f=g+1$ and $c>500000$, then (2.5) has no solutions $(x,y,z)$ with $y>z>x$.

\subsection{A shuffle variant of the Je\'smanowicz conjecture}  %2.3  

Let $(A_1, B_1, C_1)$ be a primitive Pythagorean triple with $B_1 \equiv 0 \bmod 2$.  Then $A_1$, $B_1$, and $C_1$ can be expressed as in (2.3).  In 2011, T. Miyazaki \cite{Miy6} proposed the following conjecture:

\begin{Conjecture}  % 2.2
If $f = g+1$, then the equation 
$$(f^2 + g^2)^x + (2fg)^y = (f^2-g^2)^z, x,y,z \in \natN \eqno{(2.8)} $$
has only one solution $(x,y,z) = (1,1,2)$, otherwise (2.8) has no solutions $(x,y,z)$.  
\end{Conjecture}  

This is an unsolved problem as well.  Up to now, it has been verified in the following cases.

(i)  (T. Miyazaki \cite{Miy6})  $f^2+g^2 \equiv 1 \bmod 2fg$.  

(ii)  (Z. Rabai \cite{Rab1})  $f^2 + g^2 \equiv 1 \bmod d$, where $d$ is the greatest odd divisor of $2fg$.  

(iii)  (B.-L. Liu \cite{BLLiu1})  $(f,g) \equiv (0,1)$, $(0,5)$, $(1,2)$, $(2,3)$, $(3,4)$, $(4,1)$, $(4,5)$, $(5,6)$, $(6,7)$, or $((7,0) \bmod 8$.  

(iv)  (Q. Feng \cite{Fen1})  $(f,g) = (2^k, 1)$ where $k$ is a positive integer.  

(v)  (X.-H. Wang and S. Gou \cite{WaGo1})  $(f,g) = (2^k, p)$, where $k$ is a positive integer and $p$ is an odd prime.

\section{Terai-Je\'smanowicz conjecture}  % 3

Let $p$, $q$, $r$ be fixed positive integers such that 
$$ a^p + b^q = c^r, \min\{p,q,r \} > 1.  \eqno{(3.1)} $$
In 1994 N. Terai \cite{Ter1} proposed an important generalization on Je\'smanowicz' conjecture as follows:

\begin{Conjecture}  %3.1
If $a$, $b$, $c$, $p$, $q$, $r$ satisfy (3.1), then (1.1) has only one solution $(x,y,z) = (p,q,r)$.  
\end{Conjecture}

In 1999 Z.-F. Cao \cite{Cao7} showed that Conjecture 3.1 is clearly false.  He suggested that the condition $\max\{a,b,c\}>7$ should be added to the hypotheses of Conjecture 3.1, and used the term \lq\lq Terai-Je\'smanowicz conjecture'' for the resulting statement.  In the same year, N. Terai \cite{Ter4} gave a similar modification for Conjecture 3.1.  However, Z.-F. Cao and X.-L. Dong \cite{CaDo1}, and M.-H. Le \cite{Le10}, \cite{Le11} independently found infinitely many counterexamples to the 
Terai-Je\'smanowicz conjecture.  Therefore, they stated the following conjecture:

\begin{Conjecture}  % 3.2 
  If $a$, $b$, $c$, $p$, $q$, $r$ satisfy (3.1), then the only solution to (1.1) with $\min\{x,y,z \}>1$ is $(x,y,z) = (p,q,r)$. 
 \end{Conjecture}

A history of the various versions of this conjecture is also given in \cite{SDCT1}  along with a history of the Je\'smanowicz' conjecture itself.  We clarify a potentially confusing point in the history of the various versions of Conjecture 3.2 as given in \cite{SDCT1}:  there it is stated that \cite{Le10}  corrects \cite{CaDo1}, whereas \cite{Le10} actually corrects \cite{Cao7}; \cite{CaDo1} also corrects \cite{Cao7}, and gives a version of the conjecture identical to that given in \cite{Le10}, which is Conjecture 3.2; in \cite{SDCT1}, Conjecture 3.2 is accidentally given twice (as Conjectures 4.3 and 4.4).    

The existing results concerning (1.1) all apply to Conjecture 3.2.  But, in general, Conjecture 3.2 is not yet solved.  In 2015, by using the upper bound (1.2) with some elementary methods, Y.-Z. Hu and M.-H. Le \cite{HuLe1} proved that if $\max\{a,b,c\}$ is large enough and $a$, $b$, $c$ satisfy certain divisibility conditions, then (1.1) has at most one solution $(x,y,z)$ with $\min\{x,y,z \}>1$.  

Recently, R. Scott and R. Styer \cite{ScSt4} proposed a more precise conjecture as follows.

\begin{Conjecture}  %3.3
Let $a<b$ with $\gcd(a,b)=1$.  Then we have $N(a,b,c) \le 1$, except for the following cases:

(i)  $(a,b,c) = (3,5,2)$, $(x,y,z) = (1,1,3)$, $(3,1,5)$, and $(1,3,7)$.

(ii)  $(a,b,c) = (3,13,2)$, $(x,y,z) = (1,1,4)$ and $(5,1,8)$.

(iii)  $(a,b,c) = (3,13,4)$, $(x,y,z) = (1,1,2)$ and $(5,1,4)$.

(iv)  $(a,b,c) = (3,13,16)$, $(x,y,z) = (1,1,1)$ and $(5,1,2)$. 

(v)  $(a,b,c) = (3,13,2200)$, $(x,y,z) = (1,3,1)$ and $(7,1,1)$.

(vi)  $(a,b,c) = (2,3,11)$, $(x,y,z) = (1,2,1)$ and $(3,1,1)$.

(vii)  $(a,b,c) = (2,3,35)$, $(x,y,z) = (3,3,1)$ and $(5,1,1)$.

(viii)  $(a,b,c) = (2,3,259)$, $(x,y,z) = (4,5,1)$ and $(8,1,1)$.

(ix)  $(a,b,c) = (2,5,3)$, $(x,y,z) = (2,1,2)$ and $(1,2,3)$.

(x)  $(a,b,c) = (2,5,133)$, $(x,y,z) = (3,3,1)$ and $(7,1,1)$.

(xi)  $(a,b,c) = (2,7,3)$, $(x,y,z) = (1,1,2)$ and $(5,2,4)$.

(xii)  $(a,b,c) = (2,89, 91)$, $(x,y,z) = (1,1,1)$ and $(13,1,2)$.

(xiii)  $(a,b,c) = (2,91,8283)$, $(x,y,z) = (1,2,1)$ and $(13,1,1)$.

(xiv)  $(a,b,c) = (3,4,259)$, $(x,y,z) = (1,4,1)$ and $(5,2,1)$.

(xv)  $(a,b,c) = (3,10,13)$, $(x,y,z) = (1,1,1)$ and $(7,1,3)$.

(xvi)  $(a,b,c) = (3,16,259)$, $(x,y,z) = (1,2,1)$ and $(5,1,1)$.

(xvii)  $(a,b,c) = (2,2^k-1,2^k+1)$, $(x,y,z) = (1,1,1)$ and $(k+2,2,2)$ where $k$ is a positive integer with $k \ge 2$.
\end{Conjecture}

Obviously, Conjecture 3.3 contains Conjecture 3.2 and Conjecture 2.1 for $n=1$.  Undoubtedly, we can look upon Conjecture 3.3 as the ultimate aim of solving (1.1).  Although the known results concerning (1.1) show Conjecture 3.3 to be true in many cases, the complete resolution of the conjecture is a very difficult problem.

\subsection{Earlier Results} %3.1

Since the 1950s, methods that arise in elementary number theory and algebraic number theory have been used to find all solutions $(x,y,z)$ of (1.1) for the following small values of $a$, $b$, $c$ with $a<b$.  

(i)  (T. Nagell \cite{Nag1})  $\max\{a,b,c\} \le 7$.  

(ii)  (A. Makowski \cite{Mak1})  $(a,b,c) = (2,11,5)$.

(iii) (T. Hadano \cite{Had1})  $11 \le \max\{a,b,c\} \le 17$ and $(a,b,c) \ne (3,13,2)$.  

(iv)  (S. Uchiyama \cite{Uch1})  $(a,b,c) = (3,13,2)$.  

(v)  (Z.-F. Cao and C.-S. Cao \cite{CaCa1})  $\max\{a,b,c\} = 14$.  

(vi)  (J.-H. Ren and J.-K. Zhang \cite{ReZh1})  $14 \le \max\{a,b,c\} \le 16$.  

(vii)  (Q. Sun and X.-M. Zhou \cite{SuZh1})  $\max\{a,b,c\} = 19$.  

(viii)  (X.-Z. Yang \cite{Yan1}) $\max\{a,b,c\} = 23$.  

(ix)  (Z.-F. Cao \cite{Cao3}, \cite{Cao5})  $a$, $b$, $c$ distinct primes with $29 \le \max\{a,b,c\} \le 97$.  

(x)  (X.-G. Guan \cite{Gua1}, \cite{Gua2})  $a=2$, $b$ and $c$ are distinct odd primes with $100 < \max\{b,c\} < 150$. 

(xi)  (Z.-F. Cao, R.-Z. Tong, and Z.-J. Wang \cite{CTW1})  $a=2$, $b$ and $c$ are distinct odd primes with $100 < \max\{b,c\} < 200$.  

(xii)  (X.-E. Zhou \cite{Zho1})   $c=2$, $a$ and $b$ are distinct odd primes with $200 < \max\{a,b\} < 300$.  

Note that Theorem 2 of R. Scott \cite{Sc1} gives a general method for finding all solutions to (1.1) for a given $(a,b,c)$ when $c$ is odd or $c=2$.  

Using the above mentioned results (i)-(ix), Z.-F. Cao \cite{Cao6} conjectured  that if $a$, $b$, $c$ are distinct primes with $\max\{a,b,c\} >7$, then $N(a,b,c) \le 1$.  But, owing to the fact that S. Uchiyama \cite{Uch1} missed the solution $(x,y,z) = (5,1,8)$ for $(a,b,c) = (3, 13, 2)$, Cao's conjecture is false.  In 1985, under the parity classification of solutions $(x,y,z)$ of (1.1), M.-H. Le \cite{Le1} proved that if $a$, $b$, $c$ are distinct primes with $a=2$, then every class has at most one solution $(x,y,z)$.  In the same paper, he proved that if $c=2$ and $a$, $b$, are distinct odd primes with $\max\{a,b\}>13$, then $N(a,b,c) \le 1$.  Using a different method, R. Scott \cite{Sc1} obtained similar results for $a$ and $b$ not necessarily prime.  

Let $a=2$ and let $b$, $c$ be a pair of twin primes with $2+b=c$.  Q. Sun and X.-M. Zhou \cite{SuZh1} proved that if $b>3$ and $b^2+b+1$ is an odd prime with $c^{b+1} \not\equiv 1 \bmod (b^2 + b +1)$, then (1.1) has only one solution $(x,y,z) = (1,1,1)$.  Afterwards, Z.-F. Cao \cite{Cao4}, H. Yang and R.-Q. Fu \cite{YaFu1} independently solved this case.  They proved that if $a=2$ and $b$, $c$ are twin primes with $c>b>3$, then (1.1) has only one solution $(x,y,z) = (1,1,1)$.  

Let $k$ be a fixed positive integer.  H. Edgar asked how many solutions $(y,z)$ there are to the equation 
$$c^z - b^y = 2^k, y,z \in \natN  \eqno{(3.2}) $$
for fixed distinct odd primes $b$ and $c$ (see R. K. Guy \cite[Problem D9]{Guy}).    
In 1987, Z.-F. Cao and D.-Z. Wang \cite{CaWa1}, using an earlier result of Z.-F. Cao \cite{Cao2}, proved that (3.2) has at most two solutions $(y,z)$.  In 2004, R. Scott and R. Styer \cite{ScSt1} improved this to at most one solution $(y,z)$.  In this last result, we can allow $b$ composite if, instead of using the result from \cite{Cao2}, we use a more general (but non-elementary) result in which $b$ is not restricted to be prime.  Such a result has been obtained independently by three different methods: the most direct approach is that of F. Luca \cite{Luc}, who uses results of Y. Bilu, G. Hanrot, and P. M. Voutier \cite{BHV}; a second proof was given by M.-H. Le \cite{Le}, using linear forms in logarithms and several auxiliary lemmas; the result can also be obtained by considering the equation $x^n + 2^\alpha  y^n = C z^2$, handled by M. A. Bennett and C. Skinner \cite[Theorem 1.2]{BeSk} using (in the words of those authors) \lq\lq combinations of every technique we have currently available.''  This last result is often useful in treating other exponential diophantine equations.      

M. Perisastri \cite{Per1}, M. Toyoizumi \cite{Toy1}, M.-H. Le \cite{Le2}, Z.-F. Cao \cite{Cao1} and \cite{Cao2}, Y.-S. Cao \cite{YSCao1}, J.-M. Chen \cite{JMChe1} gave certain solvability conditions to (1.1) for $a=2$ or $2l$ where $l$ is an odd prime.

\subsection{Generalized primitive Pythagorean triples} %3.2

For $p=q=2$ and $r>2$, by (3.1), we have 
$$ a^2 + b^2 = c^r.  \eqno{(3.3)} $$
Then $(a,b,c)$ is called a generalized primitive Pythagorean triple.  Since $\gcd(a,b) = 1$ and $r>2$, we see from (3.3) that $c \equiv 1 \bmod 4$ and $a$, $b$ have opposite parity.  We may therefore assume without loss of generality that $a \equiv 0 \bmod 2$ and $b \equiv 1 \bmod 2$. By L. J. Mordell \cite[Section 15.2]{Mor1}, every generalized primitive Pythagorean triple $(a,b,c)$ can be expressed as 
$$
\begin{array}{c}
a = f \left| \sum_{i=0}^{(r-1)/2} \binom{r}{2i} f^{r-2i-1} (-g^2)^i \right|, 
  b = g \left| \sum_{i=0}^{(r-1)/2} \binom{r}{2i+1} f^{r-2i-1} (-g^2)^i \right|,    c = f^2 + g^2, \\  
f,g \in \natN, \gcd(f,g)=1, f \equiv 0 \bmod 2 
\end{array}
\eqno{(3.4)}
$$
or 
$$ 
\begin{array}{c}
a = fg \left| \sum_{i=0}^{r/2-1} \binom{r}{2i+1} f^{r-2i-2} (-g^2)^i \right|, 
  b =  \left| \sum_{i=0}^{r/2-1} \binom{r}{2i} f^{r-2i} (-g^2)^i \right|,
  c = f^2 + g^2, \\  
f,g \in \natN, \gcd(f,g)=1, fg \equiv 0 \bmod 2, 
\end{array}
 \eqno{(3.5)} $$
according as $r \equiv 1 \bmod 2$ or not.  

Many recent works on (1.1) concern generalized primitive Pythagorean triples $(a,b,c)$.  These results can be divided into the following:

Case I.  If $r \in \{ 3,5,7,9 \}$ and $f$, $g$ satisfy one of the following conditions, then (1.1) has only one solution $(x,y,z) = (2,2,r)$.  

(i)  (N. Terai \cite{Ter1}; M.-H. Le \cite{Le5})  $r=3$, $g=1$, and $f$ satisfies certain congruence and divisibility conditions.   

(ii)  (M.-H. Le \cite{Le9})  $r=3$, $g=1$, $f \equiv 2 \bmod 4$, and $f \ge 206$.  

(iii)  (Z.-F. Cao and X.-L. Dong \cite{CaDo2}) $r=3$ and $g=1$.  

(iv)  (M.-H. Le \cite{Le15}; Y.-Z. Hu and P.-Z. Yuan \cite{HuYu2})  $r=3$ and $b$ is an odd prime power.  

(v)  (N. Terai \cite{Ter2}; X.-L. Dong and Z.-F. Cao \cite{DoCa1})  $r=5$, $g=1$, and $b$ is an odd prime.  

(vi)  (M.-H. Le \cite{Le13})  $r=5$, $g=1$, and $f \ge 542$.  

(vii)  (Z.-F. Cao and X.-L. Dong \cite{CaDo2}; Y.-Z. Hu and P.-Z. Yuan \cite{HuYu1})  $r=5$ and $g=1$.  

(viii)  (J.-P. Chen \cite{JPChe1})  $r=7$, $g=1$, and $f \ge 2.4 \times 10^9$.  

(ix)  (X.-G. Guan \cite{Gua3})  $r=7$, $g=1$, and $b$ is an odd prime.   $r=9$, $g=1$, and $f>2863$.

Case II.  If $g=1$, $a$, $b$, $f$, $r$ satisfy one of the following conditions, then (1.1) has only one solution $(x,y,z) = (2,2,r)$.  

(i)  (N. Terai \cite{Ter3}; N. Terai and K. Takakuwa \cite{TeTa1})  $r$ is an odd prime, $a \equiv 2 \bmod 4$, and $b$ is an odd prime with $b \equiv 3 \bmod 4$.  

(ii)  (Z.-F. Cao and X.-L. Dong \cite{CaDo1})  $r \equiv 1 \bmod 2$, $a \equiv 2 \bmod 4$, and $b$ is an odd prime with $b \equiv 3 \bmod 4$.  

(iii)  (Z.-F. Cao \cite{Cao7}) $r  \equiv 1 \bmod 2$, $a \equiv 2 \bmod 4$, $b \equiv 3 \bmod 4$, and $c$ is an odd prime.  

(iv)  (M.-H. Le \cite{Le14})  $r \equiv 1 \bmod 2$, $a \equiv 2 \bmod 4$, $f > 41 r^{3/2}$, and $c$ is an odd prime.  

(v)  (Z.-F. Cao and X.-L. Dong \cite{CaDo2})  $r \equiv 1 \bmod 2$, $a \equiv 2 \bmod 4$, $f \ge 200$, and $f/\sqrt{825 \log(f^2+1)-1} > r$.  

(vi)  (J. Y. Xia and P.-Z. Yuan  \cite{XiYu1})  $r \equiv 1 \bmod 2$, $a \equiv 2 \bmod 4$, $b \equiv 3 \bmod 4$, and $f / \allowbreak \sqrt{1.5 \log_3(f^2+1)-1} \allowbreak > r$. 
  
(vii)  (M. Cipu and M. Mignotte \cite{CiMi1})   $r \equiv 1 \bmod 2$, $a \equiv 2 \bmod 4$, and $b \equiv 3 \bmod 4$,
except for finitely many tetrads $(a,b,c,r)$ with $r<770$.  

(viii)  (M.-H. Le \cite{Le16})  $r \equiv 5 \bmod 8$ and 
$$ 
f > 
\begin{cases}
 r^2,\ {\rm if \ } \ r<11500, \\
\frac{2r}{\pi}, \ {\rm otherwise.} 
\end{cases}
$$

(ix)  (B. He, A. Togb\'e, and S.-C. Yang \cite{HTY1})  $r$ is an odd prime with $r \equiv 5 \bmod 8$ or $r \equiv 19 \bmod 24$, $f > 2r/\pi$.  

(x)  (M.-H. Le \cite{Le18}) $r\equiv 1 \bmod 2$ and $f> 10^6 r^6$.

(xi)  (T. Miyazaki \cite{Miy2}, \cite{Miy5})  $r \equiv 4$ or $6 \bmod 8$, $f^2 \log(2)/ \log(f^2+1) \ge r^3$.  

(xii)  (X.-G. Guan \cite{Gua4})  $r \equiv 0 \bmod 2$ and $ f > 48400 r^2 (\log(r))^2$.  

(xiii)  (M.-H. Le, A. Togb\'e, and H.-L. Zhu \cite{LTZ1})  $f> \max\{10^{15}, 2 r^3\}$.  

In 2012, F. Luca \cite{Luc1} basically solved this case.  He proved that if $(a,b,c)$ is a generalized primitive Pythagorean triple with $g=1$, then (1.1) has only one solution $(x,y,z) = (2,2,r)$, except for finitely many pairs $(r,f)$.  Two years later, T. Miyazaki \cite{Miy9} showed that the exceptional cases $(r,f)$ of Luca's result satisfy $r < 10^{74}$ and $f<10^{34}$.  

Similarly, T. Miyazaki and F. Luca \cite{MiLu1} proved that if $p \equiv 2 \bmod 4$, $q = r = 2$, and $a=f^2-1$ in (3.1), where $f$ is a fixed positive integer with $f \equiv 0 \bmod 2$, then (1.1) has only one solution $(x,y,z) = (p,2,2)$, except for finitely many pairs $(p,f)$.  

Case III.  If $r \equiv 1 \bmod 2$, and if $a$, $b$, and $c$ satisfy one of the following conditions, then (1.1) has only one solution $(x,y,z) = (2,2,r)$.  

(i)  (M.-H. Le \cite{Le8})  $a \equiv 2 \bmod 4$, $b \equiv 3 \bmod 8$, and $c$ is an odd prime power.  

(ii)  (N. Terai \cite{Ter4})  $a \equiv 2 \bmod 4$, $b \equiv 3 \bmod 8$, $b \ge 41 a$, and $\Jacobi{a}{b} = -1$ where $\Jacobi{*}{*}$ is the Jacobi symbol.  

(iii)  (N. Terai \cite{Ter4})  $a \equiv 2 \bmod 4$, $b \equiv 3 \bmod 8$, $b \ge 30 a$, and $b$ has an odd prime divisor $l$ with $\Jacobi{a}{l} = -1$.

(iv)  (Q. Feng and D. Han \cite{FeHa1})  $r \equiv \pm 3 \bmod 8$, $f > 2r/\pi$, and $b$ is an odd prime.  

(v)  (Z.-F. Cao and X.-L. Dong \cite{CaDo1})  $a \equiv 2 \bmod 4$, $b \equiv 3 \bmod 4$, and $b \ge 25.1 a$.  

(vi)  (M.-H. Le \cite{Le10})  $a \equiv 2 \bmod 4$, $b \equiv 3 \bmod 4$, and $\frac{b}{a} > 1/\sqrt{ e^{r/1856} -1}$.  This implies that $a \equiv 2 \bmod 4$, $b \equiv 3 \bmod 4$, and either $r \ge 1287$ or $b \ge 25 a$.  

(vii)  (M.-H. Le \cite{Le13})  $a \equiv 2 \bmod 4$, $b \equiv 3 \bmod 4$, $r>7200$, and $c > 3 \times 10^{37}$.  

(viii)  (J.-Y. Xia and P.-Z. Yuan \cite{XiYu1})  $a \equiv 2 \bmod 4$, $b \equiv 3 \bmod 4$, $b>a$, and $c-1$ is not a square.

\subsection{Miscellaneous Results}  % 3.3  

In this subsection we introduce the other works on (1.1).  Let $(p,q,r) = (1,1,1)$ in (3.1).  In 2012, T. Miyazaki and 
A. Togb\'e \cite{MiTo1} proved that, if $a =2$, then (1.1) has only one solution $(x,y,z) = (1,1,1)$, except for $b=89$.  Afterwards, M. Tang and Q.-H. Yang \cite{TaYa1} discussed a general equation as follows: 
$$ (2n)^x + (bn)^y = ((b+2)n)^z, x,y, z \in \natN, \eqno{(3.6)} $$
where $b$ is an odd integer with $b>1$.  They proved that if $(x,y,z)$ is a solution of (3.6) with $(x,y,z) \ne (1,1,1)$, then either $x>z \ge y$ or $y > z > x$.  Obviously, for any $b$, (3.6) has solutions $(x,y,z,n) = (1,1,1,t)$ and $(3,2,2, (b+1)/2)$, where $t$ is an arbitrary positive integer.  In 2014, Y.-H. Yu and Z.-P. Li \cite{YuLi1} proved that, for any $b$, (3.6) has only finitely many solutions $(x,y,z,n)$ with $x>z=y$.  In the same year, W.-J. Guan and S. Che \cite{GuCh1} found all solutions $(x,y,z,n)$ of (3.6) with $x>z =y$ for $b \not\equiv 7 \bmod 8$.  They proved that if $b \not\equiv 7 \bmod 8$, then (3.6) with $x > z=y$ has only the solutions
$$
(x,y,z,n) = 
\begin{cases}
\left(3,2,2,\frac{b+1}{2} \right) {\rm \ if \ } b \equiv 1 \bmod 4, \\
 \left(3,2,2,\frac{b+1}{2}\right) {\rm \ and \ } \left(5,4,4, \frac{b+1}{4} (b^2+2b+2)\right)   {\rm \ if \ } b \equiv 3 \bmod 4  {\rm \ and \ } \frac{b+1}{4} {\rm \ is\ not\ a\ square}, \\
 \left(3,2,2,\frac{b+1}{2}\right), \left(4,2,2, \sqrt{\frac{b+1}{4}}\right) , {\rm \ and \ } \left(5,4,4, \frac{b+1}{4} (b^2+2b+2)\right)   {\rm \ otherwise. }  
\end{cases}
$$
Moreover, T. Miyazaki, A. Togb\'e, and P.-Z. Yuan \cite{MTY1} deal with the equation 
$$a^x + b^y = (b+2)^z, x,y,z \in \natN  \eqno{(3.7)} $$
for $a \ne 2$ and $\gcd(a,b) = 1$.  They found all solutions $(x,y,z)$ of (3.7) with $b \equiv -1 \bmod a$.  Similarly, M.-H. Zhu and X.-X. Li \cite{ZhLi1} proved that if $a = 4$ and $c = b+4$, then (1.1) has only one solution $(x,y,z) = (1,1,1)$.  

Since $( p,q,r ) = ( 1,1,1 )$ in ( 3.1 ), we have $a+b = c$. Hence, for any positive integer $n$, the equation 
$$ ( an )^x + ( bn )^y = ( cn )^z , x,y, z \in \natN  \eqno{(3.8)} $$
has the solution $( x,y,z ) = ( 1,1,1 )$. A solution $( x,y,z )$ of (3.8) with $( x,y,z ) \ne ( 1,1,1 )$ is called exceptional. Very recently, C.-F. Sun and M. Tang \cite{SuTa1} showed that if $\min\{a,b \}>2$ and 
$( x,y,z )$ is an exceptional solution of (3.8), then either $x>z>y$ or $y>z>x$. In the same paper, they proved that if $( a,b ) = ( 3,5 )$, $( 5,8 )$, $(8,13 )$, or $( 13,21 )$, then (3.8) has no exceptional solutions. In this respect, P.-Z. Yuan and Q. Han \cite{YuHa1} proposed the following conjecture:

\begin{Conjecture} %3.4
If $a + b = c$ and $\min\{ a,b \} \ge 4$, then (3.8) has no exceptional solutions $( x,y,z )$.
\end{Conjecture} 

In \cite{YuHa1}, they proved that if $a$ and $b$ are squares with $b \equiv 4 \bmod 8$, then (3.8) has no solutions $(x,y,z)$ with $y>z>x$, in particular, if $b=4$, then Conjecture 3.4 is true.  Very recently, M.-H. Le \cite{Le20} proved that if $a$ and $b$ are squares with $a> 64 b^3$, then (3.8) has no solutions $(x,y,z)$ with $x > z > y$, in particular, if $b \equiv 4 \bmod 8$, then Conjecture 3.4 is true.  

Let $(p,q,r) = (2,k,1)$ in (3.1), where $k$ is an odd positive integer.  If $a \equiv -1 \bmod b^{k+1}$, $b \equiv c \equiv 1 \bmod 2$, and $b$, $c$, $k$ satisfy one of the following conditions, then (1.1) has only one solution $(x,y,z)=(2,k,1)$.

(i)  (N. Terai \cite{Ter5})  $k \in \{ 1,3\}$ and $b$ is an odd prime with $b \equiv 3 \bmod 4$.  

(ii)  (M.-H. Le \cite{Le12})  $b \equiv 3 \bmod 4$.  

(iii)  (B.-L. Liu \cite{BLLiu2})  $b \equiv 5 \bmod 24$.  

In 2006, M.-Y. Lin \cite{Lin1} completely solved this case.  He proved that if $(p,q,r) = (2,k,1)$, $a \equiv -1 \bmod b^{k+1}$, and $b \equiv c \equiv k \equiv 1 \bmod 2$, then (1.1) has only one solution $(x,y,z) = (2,k,1)$.  

Similarly, Z.-W. Liu \cite{ZWLiu1} proved that if $(p,q,r) = (1,k,2)$ where $k$ is an odd positive integer, $b \equiv 5 \bmod 12$, and $c$ is an odd prime with $c \equiv -1 \bmod b^{k+1}$, then (1.1) has only one solution $(x,y,z) = (1,k,2)$.  

Let $(a,b,c) = (f,f+1,f+2)$, where $f$ is a fixed positive integer with $f>1$.  By the result of W. Sierpi\'nski \cite{Sie1}, if $f = 3$, then (1.1) has only one solution $(x,y,z) = (2,2,2)$.  In this respect, B. Leszczy\'nski \cite{Les1} and A. Makowski \cite{Mak2} gave some conditions for (1.1) to have solutions. In 2009, B. He and A. Togb\'e \cite{HeTo1} completely solved this case.  They proved that if $f>3$, then (1.1) has no solutions $(x,y,z)$. 

Let $(a,b,c) = (f, 3 f^2-1, 4f^2-1)$, where $f$ is a fixed positive integer with $f>1$. In 2007, Y.-Z. Hu \cite{Hu1} proved that (1.1) has only the solution $(x,y,z) = (2,1,1)$ with $x \equiv 0 \bmod 2$.  Afterwards, L.-M. Chen \cite{LMChe1}, and B. He and A. Togb\'e \cite{HeTo2} independently showed that $(x,y,z) = (2,1,1)$ is the unique solution of (1.1).  

Similarly, Y.-Z. Hu \cite{Hu2} proved that if $(a,b,c) = (8 f^3 - 3f, 3f^2-1, 4f^2-1)$, where $f$ is a fixed positive integer with $f>3$, then (1.1) has only the solution $(x,y,z) = (2,1,3)$ with $ x \equiv 0 \bmod 2$.  B. He and S.-C. Yang \cite{HeYa1} showed that the solution $(x,y,z) = (2,1,3)$ is unique.  Moreover, S.-C. Yang and B. He \cite{YaHe1} proved that if $(a,b,c) = (8f^3+3f, 3f^2+1, 4f^2+1)$, where $f$ is a fixed positive integer, then (1.1) has only one solution $(x,y,z)=(2,1,3)$.  

Let $(a,b,c) = (2f, 2fg-1, 2fg+1)$, where $f$, $g$ are fixed positive integers with $\min\{f,g \}>1$.  In 2012, T. Miyazaki and A. Togb\'e \cite{MiTo1} proved that (1.1) has solutions $(x,y,z)$ if and only if either $f=2g$ or $g=f^2$.  Further, (1.1) has only one solution $(x,y,z)= (2,2,2)$ if $f=2g$ or $(3,2,2)$ if $g=f^2$.  

Let $(a,b,c) = ( fn^2+1, (g^2-f)n^2-1, gn)$, where $f$, $g$, $n$ are fixed positive integers.  If $f$, $g$, and $n$ satisfy one of the following conditions, then (1.1) has only one solution $(x,y,z) = (1,1,2)$.  

(i)  (N. Terai \cite{Ter7})  $(f,g)=(4,3)$, $n \le 20$ or $n \not\equiv 3 \bmod 6$.  

(ii)  (J.-P. Wang, T.-T. Wang, and W.-P. Zhang \cite{WWZ1})  $(f,g)=(4,3)$ and $n \not\equiv 0 \bmod 3$.  

(iii)  (L.-J. Su and X.-X. Li \cite{SuLi1})  $(f,g)=(4,3)$, $n>90$, and $n \equiv 3 \bmod 6$.  

(iv)  (C. Bert\'ok \cite{Ber1})  $(f,g) = (4,5)$, $20 < n \le 90$, and $n \equiv 3 \bmod 6$.

(v)  (N. Terai and T. Hibino \cite{TeHi1})  $(f,g)=(12,5)$, $n \not\equiv 17$ or $33 \bmod 40$.  

(vi)  (N. Terai and T. Hibino \cite{TeHi2})  $(f,g)=(3 l,l)$, where $l$ is an odd prime with $l < 3784$ and $l \equiv 1 \bmod 4$, $n \not\equiv 0 \bmod 3$, and $n \equiv 1 \bmod 4$.  

(vii)  (T. Miyazaki and N. Terai \cite{MiTe1})  $f = 1$, $g \equiv \pm 3 \bmod 8$, $n \equiv \pm 1 \bmod g$, and $(g,n) \ne (3,1)$.  

(viii)  (X.-W. Pan \cite{Pan1}; R.-Q. Fu and H. Yang \cite{FuYa1})  $f \equiv 4$ or $5 \bmod 8$, $g \equiv 1 \bmod 2$, $n > 6 g^2 \log g$, and $\Jacobi{f+1}{g} = -1$ where $\Jacobi{*}{*}$ is the Jacobi symbol.

Let $F = \{ F_m \}_{m=0}^\infty$ be the Fibonacci sequence, so we have 
$$F_0 = 0, F_1 = 1, F_{m+2} = F_{m+1} + F_m, m \ge 0.  $$
For any positive integer $t$ with $t \ge 2$, if the sequence $F^{(t)} = \{ F_m^{(t)} \}$ satisfies 
$$ F^{(t)}_0 = \dots = F^{(t)}_{t-1} = 0, F^{(t)}_t = 1, F^{(t)}_{m+t+1} = F^{(t)}_{m+t} + \dots + F^{(t)}_m, m \ge 0, $$
then $F^{(t)}$ is called the $t$-generalized Fibonacci sequence.  Correspondingly, $F_m$ and $F^{(t)}_m$ are called Fibonacci numbers and $t$-generalized Fibonacci numbers, respectively. In 2002, N. Terai \cite{Ter6} asked if (1.1) has only one solution $(x,y,z) = (2,2,1)$ for $(a,b,c) = (F_k, F_{k+1}, F_{2k+1})$, where $k$ is a positive integer with $k \ge 3$.  Thirteen years later, this question was answered by T. Miyazaki \cite{Miy11}.  In addition, he proved the further result that if $(a,b,c) = (F_k, F_{2k+2}, F_{k+2})$, where $k$ is a positive integer with $k \ge 3$, then (1.1) has only one solution $(x,y,z) = (2,1,2)$.  Moreover, F. Luca and R. Oyono \cite{LuOy1}, A. P. Chaves and D. Marques \cite{ChMa1} and \cite{ChMa2}, C. A. G\'omez Ruiz and F. Luca \cite{GoLu1}, N. Hirata-Kohno and F. Luca \cite{HiLu1} deal with (1.1) for some special Fibonacci numbers and generalized Fibonacci numbers $a$, $b$, and $c$.  

\subsection{A shuffle variant of the Terai-Je\'smanowicz conjecture}  % 3.4

Finally, we introduce a shuffle variant of the Terai-Je\'smanowicz conjecture proposed by T. Miyazaki \cite{Miy8} as follows:

\begin{Conjecture}  %3.5
Let $a$, $b$, $c$, $p$, $q$, $r$ be fixed positive integers satisfying (3.1).  Further, let $a<b$, $(a,b,c) \ne (2,7,3)$ and $(a,b,c) \ne (2,2^k-1, 2^k+1)$, where $k$ is a positive integer.  If $q = r = 2$ and $b+1=c$, then the equation
$$ c^X + b^Y = a^Z, X, Y, Z \in \natN  \eqno{(3.8)}  $$
has only one solution $(X,Y,Z) = (1,1,p)$.  Otherwise, (3.8) has no solutions $(X,Y,Z)$.
\end{Conjecture}    

In the same paper, T, Miyazaki proved that Conjecture 3.5 is true if $q=r=2$ and $b+1=c$.  In general, this problem has not been solved yet.

\bigskip

\noindent
Acknowledgment:  The authors thank Prof. P.-Z. Yuan for his valuable suggestions and recent papers.

\bigskip 

Maohua Le

Institute of Mathematics

Lingnan Normal College

Zhangjiang, Guangdong  524048

China

\bigskip

Reese Scott

Somerville MA USA

\bigskip 

Robert Styer

Department of Mathematics and Statistics

Villanova University

800 Lancaster Avenue

Villanova, PA  19085  USA

\end{document}